\documentclass[12pt]{amsart}
\usepackage{amsmath}

\theoremstyle{plain}
\newtheorem{thmintro}{Theorem}
\newtheorem{thm}{Theorem}[section]
\newtheorem{lem}[thm]{Lemma}
\newtheorem{lemma}[thm]{Lemma}
\newtheorem{cor}[thm]{Corollary}
\newtheorem{prop}[thm]{Proposition}

\theoremstyle{definition}

\theoremstyle{remark}
\newtheorem{rmk}[thm]{Remark}
\newtheorem{eg}[thm]{Example}
\newtheorem{conj}[thm]{Conjecture}

\numberwithin{equation}{thm}

\DeclareMathOperator{\Gal}{Gal}

\DeclareMathOperator{\Spec}{Spec}

\newcommand{\kb}{{\bar k}}
\newcommand{\lb}{{\bar \ell}}
\newcommand{\bP}{{\mathbb P}}
\newcommand{\bZ}{{\mathbb Z}}
\newcommand{\lra}{\longrightarrow}
\newcommand{\tC}{\tilde C}
\newcommand{\F}{\mathbb F}
\newcommand{\Fq}{\F_q}
\newcommand{\Fqk}{\F_{q^k}}

\newcommand{\Fqm}{\F_{q^m}}
\newcommand{\xyx}{X \times_Y X}
\newcommand{\xkx}{X \times_k X}

\newcommand{\myS}{C_1 \times_k C_2}

\newcommand{\size}{\#C(\Fq)}

\newcommand{\cU}{\mathcal {U}}
\newcommand{\cV}{\mathcal {V}}

\newcommand{\cT}{\mathcal {T}}
\newcommand{\cS}{\mathcal {S}}
\newcommand{\ta}{\tilde{a}}
\newcommand{\ts}{\tilde{\sigma}}

\newcommand{\cP}{\mathcal{P}}
\newcommand{\cQ}{\mathcal{Q}}

\newcommand{\cI}{\mathcal{I}}
\newcommand{\cJ}{\mathcal{J}}

\begin{document}
\title{Exceptional covers and bijections on rational points}
\author{Robert M. Guralnick}
\address{
  Department of Mathematics,
  University of Southern California,
  Los Angeles, CA 90089--1113
}
\email{guralnic@math.usc.edu}

\author{Thomas J. Tucker}
\address{
  Department of Mathematics,
  Hylan Building,
  University of Rochester,
  Rochester, NY 14627
}
\email{ttucker@math.rochester.edu}

\author{Michael E. Zieve}
\address{
Center for Communications Research,
805 Bunn Drive,
Princeton, NJ 08540
}
\email{zieve@math.rutgers.edu}

\thanks{The first author was partially supported by
NSF grant DMS-0140578.}


\begin{abstract}
  We show that if $f\colon X \lra Y$ is a finite, separable morphism
  of smooth curves defined over a finite field $\Fq$, where $q$ is
  larger than an explicit constant depending only on the degree of $f$
  and the genus of $X$, then $f$ maps $X(\Fq)$ surjectively onto $Y(\Fq)$
  if and only if $f$ maps $X(\Fq)$ injectively into $Y(\Fq)$.  Surprisingly,
  the bounds on $q$ for these two implications have different orders of
  magnitude.  The main tools used
  in our proof are the Chebotarev density theorem for covers of curves
  over finite fields, the Castelnuovo genus inequality, and ideas from
  Galois theory.
\end{abstract} 

 \maketitle


\section{Introduction}

Let $X$ and $Y$ be normal, geometrically irreducible varieties over
$\Fq$, and let $f\colon X \lra Y$ be a finite,
generically \'etale $\Fq$-morphism.
Then $f$ is called an {\em exceptional cover} if the diagonal is the
only geometrically irreducible component of the fiber product $\xyx$
which is defined over $\Fq$.  

The prototypical examples of exceptional covers are isogenies of
abelian varieties, which are exceptional whenever zero is the only
$\Fq$-rational point in the kernel.  Other families of examples will
be discussed in Section \ref{examples}.

The primary interest of exceptional
covers is that they induce bijections on rational points:

\begin{thmintro}
\label{bij}
If $f$ is exceptional, then $f$ maps $X(\Fq)$ bijectively onto $Y(\Fq)$.
\end{thmintro}

\noindent
This result is due to Lenstra (unpublished).  Special cases and weaker
versions were previously proved by Davenport and Lewis~\cite{DL},
MacCluer~\cite{Ma}, Williams~\cite{Wi}, Cohen~\cite{Co}, and
Fried~\cite{Fr,Fr3,FGS}.  (See \cite{LMZ} for variants of this result
over infinite constant fields.)

Note that, if $f$ is exceptional over $\Fq$, then $f$ is also
exceptional over $\Fqm$ for infinitely many $m$.  Thus, $f$ induces a
bijection $X(\Fqm)\to Y(\Fqm)$ for infinitely many $m$.  This unusual
property is the most important feature of exceptional covers.

In the present paper we show that this property characterizes exceptional
covers.  More precisely, we show (in Prop.\ \ref{end}) that $f$ is
exceptional if $X(\Fqm)\to Y(\Fqm)$ is either injective or surjective for
a single sufficiently large $m$.  We can make this completely explicit in case
$\dim X=1$, where it suffices to test a single $m$ larger than an
explicit constant depending only on $q$, the genus of $X$, and the
degree of $f$:

\begin{thmintro}\label{the main}
\label{main}
Let $X$ be a curve of genus $g_X$, and let $n$ be the degree of $f$.
\begin{enumerate}
\item Suppose $f$ maps $X(\Fq)$ injectively into $Y(\Fq)$, and
$\sqrt{q}>2n^2+4ng_X$.
Then $f$ is exceptional, and therefore bijective on rational points.
\item Suppose $f$ maps $X(\Fq)$ surjectively onto $Y(\Fq)$, and
  $\sqrt{q} > n!(3g_X + 3n)$.  Then $f$ is exceptional, and therefore
bijective on rational points.
\end{enumerate}
\end{thmintro}

Note that the bound in (1) is quite different from the bound in (2):
the bound in (1) is a degree-2 polynomial in $n$, while the bound in
(2) depends on $n!$.  The reason we get such a better bound under the
injectivity assumption is that injectivity is equivalent to the
nonexistence of non-diagonal rational points on components of
$X\times_Y X$, and these components have genus less than $n^2 + 2 g_X
n$.  There seems to be no curve playing an analogous role for
surjectivity, so we are forced to work on the Galois closure of
$\Fq(X)/\Fq(Y)$, which may have genus on the order of $n!(g_X+1)$.  We
do not know whether this phenomenon is indicative of the true situation,
or merely an artifact of our proof.  It is possible that there would
be counterexamples to (2) if we replaced $n!$ by any
polynomial in $\bZ[n]$.  However, we do not know any examples of
non-exceptional maps $f$ which are surjective on $\Fq$-points with
$\sqrt{q}>2n^2+4ng_X$.

Our proof of (1) uses the Weil lower bound on the number of rational
points on a curve and Castelnuovo's bound on the arithmetic genus of
curves in $X \times X$ to show that there are nondiagonal rational
points in $X \times_Y X$ when $\sqrt{q}>2n^2+4ng_X$.  Our proof of (2)
analyzes the decomposition and inertia groups of places of the Galois
closure of $\Fq(X)/\Fq(Y)$, using an analog of Chebotarev's density
theorem to translate injectivity, surjectivity, and exceptionality
into group-theoretic properties which are shown to be equivalent via
purely group-theoretic arguments.  This proof shows that, if
$\sqrt{q}>n!(3g_X+3n)$, then surjectivity and injectivity of $f$ are
equivalent to one another and to exceptionality.  By contrast, our
proof of (1) does not directly yield surjectivity of $f$ (although
surjectivity follows by combining (1) with Theorem \ref{bij}).

Results along the lines of Theorem~\ref{main} were previously proved
in the case $g_X=0$, where of course injectivity and surjectivity are
equivalent.  Most previous work restricts further to the case where $g_X=0$
and some point of $Y(\Fq)$ is totally ramified under $f$.  In this
case a noneffective version of our result was proved by Davenport and
Lewis~\cite{DL}.
  The best known effective version says that, if $f$
is bijective on rational points and $q\geq n^4$, then $f$ is
exceptional~\cite[pp.~51--52]{FJ} (see \cite{LY} for corrections to
\cite{FJ}).  When $g_X=0$ our result draws this conclusion under the
assumption $q\geq 4n^4$; but it is easy to modify our argument to make
use of the ramification assumption and recover the $q\geq n^4$ bound.
The effective version of the Davenport-Lewis argument extends at once
to the general case $g_X=0$ (no longer assuming a totally ramified
rational point), giving the bound $q\geq 16 n^4$~\cite{GM}.  Our
result improves this to $q\geq 4n^4$.

We prove (1) in the next section, using an argument which is similar
in spirit to that of Davenport and Lewis, although with several new
ingredients to address difficulties new to the case $g_X>0$.  In
Section~\ref{galois theory} we prove (2) and Theorem~\ref{bij}, using
a Galois-theoretic setup we learned from Lenstra.  Our proof of (2)
uses an analog of Chebotarev's density theorem, which we prove in
Section~\ref{L-C}.  We conclude in Section~\ref{examples} with some
examples and conjectures.

Let us say a few words about the terminology in this paper.  Given a
variety $W$ over a field $k$ and an extension $k'$ of $k$, we let
$W(k')$ denote the set of $k$-morphisms from $\Spec k'$ into $W$.
In particular, $W(k)$ is the set of closed points of $W$ with
residue field $k$.  Also $\bar k$ denotes an algebraic closure of $k$.
Throughout this paper, all curves are assumed to be projective and
geometrically integral.  Not all curves are assumed to be smooth; some
of the curves we work with in Section~\ref{geometry} may be singular.
\\
\\
\noindent {\it Acknowledgments.}  We would like to thank Hendrik
Lenstra for his generous help.  In particular, his ideas permeate
Sections \ref{L-C} and \ref{galois theory}.


\section{Geometry}\label{geometry}

In this section we use a geometric approach.  Our main result concerns
maps and curves defined over the field $\Fq$.  The first few
propositions are valid over an arbitrary ground field $k$ and are
stated as such.  Throughout this section, $f\colon X \lra Y$ is a finite separable
morphism between smooth curves, and $f$ is defined over either $k$ or
$\Fq$ depending on context.  We denote the geometric genus of a curve
$C$ as $g_C$ and the arithmetic genus as $p_a(C)$.  Finally, by
`component' we always mean geometric component.

Our first result shows that Castelnuovo's upper bound on the geometric
genus of a curve on a split surface is also an upper bound on the
arithmetic genus.

\begin{prop}\label{Cast}
  Let $C_1$ and $C_2$ be smooth curves and let $C'$ be a curve for which
  there is a generically injective map $\phi:C' \lra C_1 \times_k
  C_2$.  For $i=1,2$, let $g_i$ be the genus of $C_i$, let $\pi_i$
  denote projection from $\myS$ onto its $i$-th factor, and let $d_i$
  be the degree of the map $\pi_i \circ \phi: C' \lra C_i$.  Then
\begin{equation}\label{Cast eq}
p_a(\phi(C')) \leq (d_1 - 1)(d_2 - 1) + d_1 g_1 + d_2 g_2.
\end{equation}   
\end{prop}
\begin{proof}
We use several results from \cite[\S V.1]{H}, which is the source of
all references in this proof.  For divisors $D_1$ and $D_2$, denote
the intersection pairing by $D_1.D_2$.  By Thm.\ 1.1, this pairing
is symmetric, additive, and depends only on the linear equivalence
class of each $D_i$.  Let $F_i$ be a fiber of $\pi_i$.  Since $F_1$ is
linearly equivalent to any other (disjoint) fiber of $\pi_1$, we have
$F_1.F_1=0$.
The adjunction formula (Prop.\ 1.5) implies $2g_2-2=F_1.K$, where
$K$ is the canonical divisor on $\myS$.  Next, ex.\ 1.5 says that
$K.K=8(g_1-1)(g_2-1)$, so $K.K=2(F_1.K)(F_2.K)$.
By ex.\ 1.9, $K$ is numerically equivalent to $(2g_1-2)F_1+(2g_2-2)F_2$.

Let $D=\phi(C')$.  By ex.\ 1.9, $D.D\le 2d_1d_2$, so ex.\ 1.3
implies $2p_a(D)-2\le 2d_1d_2+D.K$.  Since $K\equiv
(2g_1-2)F_1+(2g_2-2)F_2$, we have $D.K=(2g_1-2)d_1+(2g_2-2)d_2$.
Thus
$$
2p_a(D)-2 \le 2d_1d_2 + (2g_1-2)d_1 + (2g_2-2)d_2,
$$
and the desired result follows.
\end{proof}

Write $Z = \xyx$, and note that $Z$ embeds naturally into $\xkx$ as
the locus of points $(P,Q)$ for which $f(P)=f(Q)$.
\begin{prop}\label{INT}
Let $(P,Q) \in Z(\kb)$ be a point which lies in more than one component
of $Z$.  Then $f$ is ramified at both $P$ and $Q$.
\end{prop}

\begin{proof}
If $f$ is unramified at $P$ then $f$ is \'etale at $P$, so $f$ is
smooth on an open subset $U$ of $X$ containing $P$.  Since the projection
$\pi_1\colon Z\lra X$ is obtained from $f$ by base extension,
it follows that $U\times_Y X$ is smooth over $X$ (by \cite[Prop. III.10.1]{H})
and thus over $k$.  This contradicts the fact that $(P,Q)$ lies
in multiple components.
\end{proof}

\begin{cor}\label{diag}
  Let $f\colon X \lra Y$ be a finite separable morphism of smooth curves and
  suppose that $f$ is injective on $X(k)$.  Then any component of
  $\xyx$ other than the diagonal contains at most $(2g_X + 2 \deg f -
  2)$ $k$-rational points.
\end{cor}
\begin{proof}
  Let $D$ be a nondiagonal component of $Z=\xyx$.  Since $f$ is injective on
  $X(k)$, every point in $Z(k)$ has the form $(P,P)$; hence
  $D(k)$ lies in the support of the intersection of $D$
  with the diagonal.  It follows that all points in $D(k)$ have the form $(P,P)$
  where $f$ ramifies at $P$.  By the Riemann-Hurwitz theorem, the number of
  such $P$ is at most
\begin{equation*}
2g_X - 2 - \deg f(2g_Y - 2) \leq 2g_X + 2 \deg f - 2.
\end{equation*}
This completes the proof.
\end{proof}

Our next result generalizes the Weil bound for the number of
$\Fq$-rational points on a smooth curve to the case of an arbitrary
curve.

\begin{prop}\label{sing-weil}
For any curve $C$ over\/ $\Fq$, we have
\begin{equation*}
\mid \#C(\Fq) - q - 1\mid \leq 2 p_a(C) \sqrt{q}.
\end{equation*}
\end{prop}

\begin{proof}
  Let $\tC$ be the normalization of $C$.  Then $\tC$ is regular (by
  \cite[Thm.\ 11.2]{Mats}) and therefore smooth (by
  \cite[Thms.\ 25.2 and 25.3]{Mats}).  The normalization map
  $\tC\lra C$ is an isomorphism away from at most $p_a(C)-g_C$
  points of $\tC$, so
\begin{equation*}
\mid \#\tC(\Fq) - \#C(\Fq)\mid \leq p_a(C) - g_C.
\end{equation*}
Since $\tC$ is smooth, we also have~\cite{We}
\begin{equation*}
\mid \#\tC(\Fq) - q - 1\mid \leq 2 g_C \sqrt{q}.
\end{equation*}
Our result follows.
\end{proof}

\begin{thm}\label{main-thm}
  Let $f\colon X \lra Y$ be a finite separable morphism of degree $n \geq 2$
  between smooth curves over\/ $\Fq$.  Suppose that $f$
  induces an injection from $X(\Fq)$ into $Y(\Fq)$ and that
\begin{equation}\label{quad}
\sqrt{q} > 2 (n-2)^2 + 4 (n-1) g_X + 1.
\end{equation}
Then $f$ is exceptional.
\end{thm}

\begin{proof}
  We argue by contradiction.  Suppose that $f$ is not exceptional.
  Then there exists a non-diagonal geometric component $C$ of $\xyx$
  that is defined over $\Fq$.  Since the projection maps from $\xyx$
  onto $X$ are generically $n$-to-1, one sees that the projection maps
  restricted to $C$ are generically at most $(n-1)$-to-1 (since the
  diagonal is also a component of $\xyx$).  Proposition
  \ref{Cast} implies that $p_a(C) \leq (n-2)^2 + 2g_X(n-1)$.
Now Proposition \ref{sing-weil} gives
$$\size \geq q + 1 - 2p_a(C)\sqrt{q} \geq q + 1 - 2 ((n-2)^2 +
2g_X(n-1)) \sqrt{q}.$$

On the other hand, by Corollary~\ref{diag} we have
$\size \leq 2 g_X + 2n - 2$.  Finally,
it is easily checked that (\ref{quad}) implies
$$q + 1 - 2 ((n-2)^2 + 2g_X(n-1)) \sqrt{q}  - (2 g_X + 2n - 2) > 0,$$
so we have our contradiction.
\end{proof}

\begin{rmk}
  The proof of Theorem~\ref{main-thm} can be modified to work under
  the weaker hypothesis that $f$ is injective over non-branch points
  of $Y$.  Injectivity is used only in Corollary~\ref{diag}; since
  this weaker hypothesis still implies that if $(P,Q)$ is a
  $k$-rational point on $\xyx$ then $f(P)=f(Q)$ is a branch point of
  $f$, we can replace the bound $(2g_X+2\deg f-2)$ from
  Cor.~\ref{diag} with the bound $(2 g_X + 2\deg f - 2)(\deg f - 1)$.
  This enlarges the bound \eqref{quad} slightly (adding $2$ to the
  right hand side is sufficient), but otherwise the reasoning is
  identical.
\end{rmk}


\section{Chebotarev}\label{L-C}

In this section we prove analogs of Chebotarev's density theorem for
normal varieties over a finite field, which we apply in the next two
sections.

Let $R$ be a commutative ring, and let $A$ be a group of automorphisms of
$R$.  We denote the fixed ring $R^A$ as $B$, and we say that $R$ is
a Galois extension of $B$.  For a single element $a \in A$ we write
$R^a$ instead of $R^{\langle a \rangle}$.  Fix a prime $\cQ$ in $R$
lying over a prime $\cP$ in $B$, and let $D=D(\cQ/\cP)$ and
$I=I(\cQ/\cP)$ denote the decomposition and inertia groups at $\cQ$.
If $\cJ$ is a prime ideal in the commutative ring $Z$, we write $m_\cJ$
for the field of fractions of $Z/\cJ$.

The following result is standard and easy
(e.g., see \cite[Thm.  2, p.~331]{Bo}).

\begin{lemma}\label{from-Bo}
  Suppose that $m_\cP$ is perfect. Then
\begin{enumerate}
\item $A$ is transitive on the set of primes of $R$ lying over $\cP$;
\item $m_\cQ/m_\cP$ is a finite Galois extension of degree
  $[D:I]$; and
\item $D/I \cong \Gal(m_\cQ/m_\cP)$.
\end{enumerate}
\end{lemma}

The next result is also known (e.g., see \cite{Wa}), but we
include a proof for the sake of completeness.  Let $H$ be a subgroup
of $A$, let $U=R^H$, and let $\cS$ be the set of left cosets of $H$ in $A$.
We may view $A$ as a group of permutations of the set $\cS$.

\begin{lemma}\label{from-Wa}
  The number of primes $\cJ \subset U$ lying over $\cP$ such that
  $m_\cJ = m_\cP$ is equal to the number of common orbits of $D$ and
  $I$ on $\cS$.  In particular, if $I = 1$ then the number of primes
  $\cJ \subset U$ lying over $\cP$ such that $m_\cJ = m_\cP$ is equal
  to the number of fixed points of $D$ on $\cS$.
\end{lemma}
\begin{proof}
For $a\in A$, let $\cQ' = a^{-1}\cQ$ and $\cJ=\cQ'\cap U$.
Since $H \cap a^{-1}Da$ and $H \cap a^{-1}Ia$ are the decomposition
and inertia groups for $\cQ'$ over $\cJ$, we have
\begin{equation*}
[m_\cJ:m_\cP] = \frac{[m_{\cQ'}:m_\cP]}{[m_{\cQ'}:m_\cJ]} =
 = \frac{[a^{-1}Da:a^{-1}Ia]}{[H \cap a^{-1}Da: H  \cap a^{-1} I a]}.  
\end{equation*}
Using the fact that $|H \cap a^{-1} M a| = |M| |H| / |MaH|$ for any
subgroup $M$ of $G$, we thus obtain
$$
[m_\cJ:m_\cP] = [D:I] \cdot \frac{|DaH|}{|D||H|} \cdot \frac{|I||H|}{|IaH|} =
\frac{|DaH|}{|IaH|},$$
which is equal to 1 if and only if $DaH$ is a common orbit of $D$ and
$I$.

For $b\in A$, we have $(b^{-1}\cQ)\cap U = \cJ$ if and only if
$DaH=DbH$.  Thus, we achieve the desired result by summing over all orbits
$DbH$ of $D$ on $\cS$.
\end{proof}

Let $W$ be a normal variety over the finite field $k$, and let $V$ be a
normal variety over a finite extension $\ell$ of $k$.  Let
$\rho: V \lra W$ be a finite, generically \'{e}tale map of $k$-schemes.
Write $K$ and $L$ for the fields of rational functions on $W$ and $V$,
so that $\rho$ induces an inclusion $K \hookrightarrow L$.
Assume that $L/K$ is Galois, and put $A=\Gal(L/K)$ and
$G=\Gal(L/K.\ell)$ (here $K.\ell$ denotes the compositum of $K$ and
$\ell$ in $L$).  Then $A/G\cong\Gal(\ell/k)$ is cyclic.
Pick $a\in A$ with $\langle aG\rangle=A/G$.  Let
$t$ be an extension of $k$ such that $[t:k]=\#\langle a\rangle$.  Note
that $t$ contains $\ell$.  Pick an automorphism $\ta$ of the
compositum $L.t$ such that $\ta|_L=a$ and $t^{\ta} =k$; such an
automorphism exists because $\ell^a=k$.  Then $(L.t)^{\ta}\supseteq
R_i^a\supseteq B_i$, and $k$ is algebraically closed in $(L.t)^{\ta}$.

Galoisness of $L/K$ implies that $V^A = W$, in the sense that
$W$ admits an affine cover $M_i = \Spec B_i$ such that $\rho^{-1}(M_i)
= \Spec R_i$ and $R_i^A = B_i$.  Then each $R_i$ is normal, so each
$B_i$ is as well (\cite[V.1.9]{Bo}).  Furthermore, $R_i$ is the
integral closure of $B_i$ in $L$ since $R_i$ is normal and integral
over $B_i$.  The ring $T_i = R_i.t$ is mapped to itself by $\ta$.
Let $V_t$ be the variety obtained from $V$ by base-extension from
$\ell$ to $t$, and let $V_t^{\ta}$ be the quotient variety of $V_t$
obtained by piecing together the fixed rings $T_i^{\ta}$.
 
The degree of the field $L.t$ over the field of fractions of
$T_i^{\ta}$ is equal to $\#\langle \ta\rangle=\#\langle
a\rangle=[t:k]$, so $T_i^{\ta}.t$ has field of fractions $L.t$.  Now,
$T_i^{\ta}.t$ and $T_i$ are both normal, because $R$ and $T_i^{\ta}$
are normal (\cite[6.7.4]{EGA}).  Both $T_i$ and $T_i^{\ta}.t$ are
integral over $T_i^{\ta}$ as well, so we must have $T_i =
T_i^{\ta}.t$.  Since $T_i$ is a finite Galois extension of both
$T_i^{\ta}$ and $B_i.t$, it is also a finite Galois extension of
$T_i^{\ta}\cap B_i.t=B_i$.

We define the degree of a maximal ideal $\cI$ in any of the rings $B$,
$R$, $T_i$, $T_i^{\ta}$ to be $[m_\cI:k]$.  Let $J \in V_t^{\ta}(k)$
and let $\cJ$ be the corresponding degree one maximal ideal in some
$T_i^{\ta}$.  Then $T_i / T_i \cJ \cong k \otimes_\ell t \cong t$, so
$T_i \cJ$ is the unique prime in $T_i$ lying over $\cJ$.  Thus, the
map $\phi_i:\cJ \mapsto T_i \cJ \cap R_i$ gives a well-defined map
from degree one maximal ideals in $T_i^{\ta}$ to maximal ideals in
$R_i$.

\begin{lem}\label{affine}
  Let $\cQ$ be a maximal ideal in $R_i$ that lies over a degree-one
  maximal ideal $\cP$ of $B_i$.  If $\langle a I(\cQ/\cP) \rangle =
  D(\cQ/\cP)/I(\cQ/\cP)$, then there are exactly $[m_\cQ:\ell]$ degree
  one maximal ideals $\cJ \in T_i^{\ta}$ such that $\phi_i(\cJ) =
  \cQ$.  Otherwise, there are no degree one maximal ideals $\cJ \in
  T_i^{\ta}$ such that $\phi_i(\cJ) = \cQ$.
\end{lem}
\begin{proof}  
  Suppose that $\langle a I(\cQ/\cP) \rangle = D(\cQ/\cP)/I(\cQ/\cP)$.
  Since $a \cQ = \cQ$, we must have $\ta T_i \cQ = T_i \cQ$.  Thus
  $\ta$ acts on $T_i / T_i \cQ \cong m_\cQ \otimes_\ell t$ by acting as
  $a$ acts on $m_\cQ$ and as $\ta$ acts on $\ell$.  The primes in
  $T_i$ lying over $\cQ$ correspond to the primes in $m_\cQ
  \otimes_\ell t$.  Now, since $a$ generates $\Gal(m_\cQ/k)$ and $\ta$
  generates $\Gal(t/k)$, there is a map $\psi:m_\cQ \hookrightarrow t$
  such that $\psi(a x) = \ta \psi(x)$.  Then the $[m_\cQ:\ell]$ primes in
  $m_\cQ \otimes_\ell$ correspond to the kernels of the maps
  $p_j:m_\cQ \otimes_\ell t \lra t$ given by $p_j(u \otimes v) =
  (\psi(a^{[\ell:k]j} u) v)$ for $0 \leq j \leq [m_\cQ:\ell] - 1$.  Since the
  kernel of $p_j$ is the set of all $\sum_n(u_n \otimes v_n)$ such
  that $\sum_n \psi(a^{[\ell:k]j} u_n) v_n = 0$, the kernel of $p_j$ is
  preserved by the action of $\ta$, so $\ta \cQ' = \cQ'$ for all
  $\cQ'$ lying over $\cQ$.  Writing $\cJ = \cQ' \cap T_i^{\ta}$, we
  then have $\langle \ta \rangle = D(\cQ'/\cJ)$ since $T_i$ is
  unramified over $T_i^{\ta}$.  For each of these $[m_\cQ:\ell]$ maximal
  ideals $\cJ$, we have $\phi_i(\cJ) = \cQ$.
  
  Now, suppose that $\langle a I(\cQ/\cP) \rangle \not=
  D(\cQ/\cP)/I(\cQ/\cP)$.  Let $\cQ'$ be a maximal ideal in $T_i$ such
  that $\cQ' \cap R_i = \cQ$ and let $\cJ = \cQ' \cap T_i^{\ta}$.  If
  $a \notin D(\cQ/\cP)$, then $\ta \cQ' \not= \ta \cQ'$, so there is
  more than one prime in $T_i$ lying over $\cJ$, which means that
  $\cJ$ cannot have degree one.  If $a \in D(\cQ/\cP)$ but does not
  generate $D(\cQ/\cP)/I(\cQ/\cP)$, then $\cQ \cap R_i^a$ has degree
  greater than one, by (3) of Lemma \ref{from-Bo}, so $\cJ$ does also,
  since $\cJ$ lies over $\cQ \cap R_i^a$.
\end{proof}
  
If $Q$ and $P$ are closed points of $V$ and $W$ with $\rho(Q) = P$, we
denote the
decomposition and inertia groups of $Q$ over $P$ as $D(Q/P)$ and
$I(Q/P)$, respectively.  Clearly, these are the same as $D(\cQ/\cP)$
and $I(\cQ/\cP)$ where $\cQ$ is a prime in some $R_i$ corresponding to
$Q$ and $\cP$ is a prime in $B_i$ such that $\cQ \cap B_i =
\cP$.  Similarly, we define $m_Q$ to be $m_\cQ$.
     
\begin{prop}\label{general cheb}
With notation as above,
\begin{equation}\label{gen}
\sum\limits_{P \in W(k)} \sum_{\substack{\rho(Q) = P \\ \langle a
    I(Q/P) \rangle = D(Q/P)}}[m_Q:\ell] =
\# V_t^{\ta}(k)
\end{equation}
where $I(Q/P)$ is the inertia group of $Q$ over $P$.
\end{prop}
\begin{proof}
  The maps $\phi_i$ patch together to form a map $\phi $ from $J(k$)
  to closed points in $V$; indeed if we let $\phi$ be the map that
  takes a point $J \in W_t^{\ta}(k)$ to the closed point $Q$ of $W$
  lying under the unique closed point of $W_t$ that lies over $J$,
  then $\phi$ agrees with $\phi_i$ on each affine piece $\Spec
  W_t^{\ta}$.  The proposition thus follows from Lemma~\ref{affine}.
\end{proof}
  
\begin{cor}\label{smooth}
  Suppose that $V$ and $W$ are nonsingular and projective.  Let $r =
  \dim V$ and let $b_0,\dots,b_{2r}$ be the Betti numbers (see
  \cite[p.~451 and 456 ]{H}) of $V$.  Then
  $$
  \left| \Big( \sum\limits_{P \in W(k)} \sum_{\substack{\rho(Q) = P
        \\ \langle a I(Q/P) \rangle =
        D(Q/P)}}[m_Q:\ell] \Big) - (\#k)^r \right| \leq \left|
    \sum\limits_{i=0}^{2r-1} (-1)^i b_i (\#k)^{i/2} \right|$$
\end{cor}
\begin{proof}
  Since $T_i^{\ta}.t = T_i$ on each affine piece $\Spec T_i$ of $V_t$,
  we see that $V_t^{\ta}$ with the base extended from $k$ to $t$ is
  isomorphic to $V_t$ (i.e., $\left( V_t^{\ta} \right)_t \cong V_t$).
  Thus, $V_t^{\ta}$ is also nonsingular (\cite[6.7.4]{EGA}), and $V$,
  $V_t$, and $V_t^{\ta}$ all have the same Betti numbers.  Thus,
  applying the Weil bound (\cite{Deligne}, see also \cite[Appendix
  3]{H} for an overview) to $V_t^{\ta}(k)$ in \eqref{gen} gives the
  desired result.
\end{proof}      
 
Proposition \ref{general cheb} also gives rise to a generalization of
the effective Chebotarev density theorem for curves that Murty and
Scherk proved in \cite{MS} (see also \cite[Chapter 5]{FJ}).  Let $\cV$
denote the set of all unramified points in $V(\lb)$ that lie over
points in $W(k)$; let $\cV_a$ denote the set of all points in $\cV$
that correspond to closed points $Q$ of $V$ such that $I(Q/\rho(Q))$
is trivial and $\langle a \rangle = D(Q/\rho(Q))$.  Note that counting
points in $V(\lb)$ is different from counting closed points; each
closed point $Q$ on $V$ corresponds to $[m_Q:\ell]$ distinct points in
$V(\lb)$.

  

\begin{cor}\label{chebo}
  Suppose that $V$ and $W$ are nonsingular and projective.  Let $r =
  \dim V$, let $b_0,\dots,b_{2r}$ be the Betti numbers of $V$, let
  $c_0,\dots,c_{2r}$ be the Betti numbers of $W$, and let $U$ be the
  ramification locus of $\rho$ in $W$ (thought of as a subscheme of
  $W$).  Then
\begin{equation}\label{our cheb}
\begin{split}
  \left| \# \cV_a - \frac{\# \cV}{\#G} \right| \leq &  
  (\#G)(\#U(k)) + \left| \sum\limits_{i=0}^{2r-1} (-1)^i b_i
    (\#k)^{i/2}
  \right| \\
   & + \left| \sum\limits_{i=0}^{2r-1} (-1)^i c_i (\#k)^{i/2} \right|.
\end{split}
\end{equation}
\end{cor}

\begin{proof}
  If $I(Q/P)$ is trivial than $\rho$ does not ramify at $Q$, so $Q$
  does not lie over a point in the ramification locus of $\rho$.  As
  noted above, each such $Q$ corresponds to $[m_Q:\ell]$ points in
  $\cV_a$.  Letting $\cU_a$ denote the set of $J \in V_t^{\ta}(k)$
  lying over points in $U(k)$ and applying Proposition \ref{general
    cheb}, we obtain $ \# \cV_a = \#V_t^{\ta}(k) - \# \cU_a,$ Since
  the degree of $V_t^{\ta}$ over $W$ is $\#G$, we have $\# \cU_a \leq
  (\#G)(\#U(k))$.  Thus, the Weil bound for $V_t^{\ta}(k)$ yields
\begin{equation}\label{V}
\begin{split}
& (\#k)^r -  \left| \sum\limits_{i=0}^{2r-1} (-1)^i b_i (\#k)^{i/2}
\right|  - (\#G)(\#U(k)) \\
& \leq \cV_a \\
  & \leq (\#k)^r + \left|
\sum\limits_{i=0}^{2r-1} (-1)^i b_i (\#k)^{i/2} \right|. 
\end{split}
\end{equation}
Similarly, we obtain
\begin{equation}\label{W}
\begin{split}
& (\#G) \left((\#k)^r -  \left| \sum\limits_{i=0}^{2r-1} (-1)^i c_i (\#k)^{i/2}
\right|  - (\#U(k)) \right) \\
& \leq \cV \\
  & \leq (\# G) \left( (\#k)^r + \left|
\sum\limits_{i=0}^{2r-1} (-1)^i c_i (\#k)^{i/2} \right| \right),
\end{split}
\end{equation}
by using the Weil bound for $W(k)$.  Dividing \eqref{W} by $\#G$ and
subtracting it from \eqref{V} yields \eqref{our cheb}.   
\end{proof}

\begin{rmk}
  When $V$ and $W$ are smooth curves, Corollary \ref{chebo} is a
  slight improvement of \cite[Theorem 1]{MS}.  Note that in this case,
  the ramification locus corresponds to a finite set of points in
  $W(\kb)$.  To make Corollary \ref{chebo} completely explicit in the
  higher-dimensional case, one must use bounds on
  $U(k)$ (such as those that come from applying the Weil bounds to
  desingularizations of the components of $U$, for example).
\end{rmk}

Recall that $g_C$ denotes the genus of a curve $C$.

\begin{cor}\label{Lenstra Cheb}
  Suppose that $W$ and $V$ are smooth projective curves.  Let $\cU$ be
  a finite subset of\/ $W(k)$, and pick $a\in A$ with $\langle
  aG\rangle=A/G$.  If
\begin{equation}\label{q bound}
\sqrt{\#k}\ge 2g_V + \sqrt{(\#G)(\#\cU)},
\end{equation}
then there is a closed point $Q$ on $V$ lying over a point $P \in
W(k)\setminus\cU$ such that $a \in D(Q/P)$ and $\langle a
I(Q/P)\rangle = D(Q/P)/I(Q/P)$.
\end{cor}

\begin{proof}
The Weil bound says that
\begin{equation*}
V_t^{\ta}(k) \ge \#k + 1 - 2g_{V_t^{\ta}}\sqrt{\#k}.
\end{equation*}
As in Corollary~\ref{smooth}, we have $g_{V_t^{\ta}} = g_V$.  Now,
(\ref{q bound}) implies that
\begin{equation*}
V_t^{\ta}(k) \ge 1 + \sqrt{\#k}\sqrt{(\#G)(\#\cU)} \geq (\#G)(\#\cU).
\end{equation*}
Since the number of closed points of $V_t^{\ta}$ that lie over points
in $\cU$ is at most $([L.t : K])(\#\cU)=(\#G)(\#\cU)$, it follows that
there is a $J \in V_t^{\ta}(k)$ that lies over a point $P \in
W(k)\setminus\cU$.
\end{proof}
   
When $V$ is singular, it is difficult to get something as uniform as
Corollary~\ref{smooth}, since we are not able to apply the Weil
bound.  It follows from the older estimate of Lang-Weil (\cite[Theorem
1]{LW}), however, that for any variety $Z$ of dimension $r$ over $k$,
there is a constant $\delta$ (depending on $Z$) such that for all
extensions $k'$ of $k$, one has
\begin{equation*}
    \left| Z(k') - (\#k')^r \right| \leq \delta (\#k')^{r - \frac{1}{2}}.
\end{equation*}

This can be proved by induction on the dimension of $Z$.  If $r = 0$,
then $Z$ is a point and we're done.  Otherwise, let $Z'$ be an affine
subset of $Z$ and let $\overline{Z'}$ be a projective closure of $Z'$.
Applying the Lang-Weil estimate to $\overline{Z'}$ and the inductive
hypothesis to $\overline{Z'} \setminus Z'$ and $Z \setminus Z'$
finishes the proof.

This allows us to treat the case of a single map $\rho:V \lra W$ with
$k$ and $\ell$ varying.  Let $k'$ be an extension of $k$ and let
$A_{K.(k'\cap \ell)}$ be the subgroup of $A$ fixing $K.(k' \cap
\ell)$.  Each element in $a \in A_{K.(k'\cap \ell)}$ extends to an
element $a_k' \in \Gal(L.k'/K.k')$ that acts as $a$ does on $L$ and
acts trivially on $k'$.  Since
$$\#\Gal(L.\ell/K.k') = \frac{\#A}{[k'\cap \ell:k]} = A_{K.(k'\cap
  \ell)},$$
every element of $\Gal(L.k'/K.k')$ is equal to $a_{k'}$
for some $a \in A_{K.(k'\cap \ell)}$.  For convenience, we denote
$\Gal(L.k'/K.k')$ as $A_{k'}$ and $\Gal(L.k'/(K.k'.\ell))$ as
$G_{k'}$.  We let $\rho_{k'}$ denote $\rho$ with its base extended to
$k'$; we have $\rho_{k'}: V_{k'.\ell} \lra W_{k'}$.

\begin{cor}\label{non smooth}
  Let $r = \dim V$.  There is a constant $\delta$ such that for any
  finite extension $k'$, we have
  $$
  \left| \Big( \sum\limits_{P \in W(k')} \sum_{\substack{\rho_{k'}(Q) \in P \\
        \langle \sigma I(Q/P) \rangle =
        D(Q/P)}}[m_Q:(\ell.k')] \Big) - (\#k')^r \right| \leq
  \delta (\#k')^{r- \frac{1}{2}}
  $$
for any $\sigma \in A_{k'}$ such that $\langle \sigma G_{'k'} \rangle
= A_{k'}/G_{k'}$. 
\end{cor}   
\begin{proof}
  Let $\sigma \in \Gal(L.k'/K.k')$.  We may write $t = a_{k'}$ for
  some $a \in A_{K.(k' \cap k)}$ as we saw above.  We define $t$ and
  $\ta$ be before and define $\ts$ to be the automorphism of $L.k'.t$
  that acts as $\sigma$ does on $L.k'$ and acts as $\ta$ does on $t$.
  Then $V_{t.k'}^{\ts}$ is isomorphic to $\left(V_t^{\ta}
  \right)_{k'}$, so, by~\eqref{gen}, we have
  \begin{equation*}
  \begin{split}
    & \left| \Big( \sum\limits_{P \in W(k')}
      \sum_{\substack{\rho_{k'}(Q) \in P \\
          \langle \sigma I(Q/P) \rangle =
          D(Q/P)}}[m_Q:(\ell.k')] \Big)  - (\#k')^r \right| \\
    & = \left|  V_{t.k'}^{\ts}(k') - (\#k')^r \right| \\
    & \leq \delta_a (\#k')^{r- \frac{1}{2}},
 \end{split}
 \end{equation*}
 for some constant $\delta_a$ depending only on $a$.  Letting $\delta$ be
 the maximum of all of the $\delta_a$ gives the desired result.
\end{proof}   
   

\section{Exceptionality}\label{galois theory}

Let $X$ and $Y$ be normal, geometrically irreducible varieties over $\Fq$,
and let $f:X \lra Y$ be a finite, generically \'{e}tale $\Fq$-morphism.
In this section we give Lenstra's Galois-theoretic proof that exceptionality
of $f$ implies bijectivity of the induced map $X(\Fq)\to Y(\Fq)$.  We then
use the same Galois-theoretic setup to show that, if $X$ and $Y$ are curves
and $q$ is sufficiently large compared to $n$ and $g_X$, then injectivity,
surjectivity, and exceptionality are equivalent.

We begin with some notation.  We will view the function field $\Fq(Y)$ as
a subfield of $\Fq(X)$, via the inclusion $\Fq(Y) \hookrightarrow \Fq(X)$
induced by $f$.  Since $f$ is
generically \'{e}tale, the extension $\Fq(X)/\Fq(Y)$ is separable.  Let
$\Omega$ be the Galois closure of this extension.  Let $\Fqk$ be
the algebraic closure of $\Fq$ in $\Omega$.  Put
$A=\Gal(\Omega/\Fq(Y))$ and $G=\Gal(\Omega/\Fqk(Y))$; then $G$ is a
normal subgroup of $A$ and $A/G\cong\Gal(\Fqk/\Fq)$ is cyclic of order
$k$.  Let $H=\Gal(\Omega/\Fq(X))$.  We may view $A$ as a group of
permutations of the set $\cS$ of left cosets of $H$ in $A$.  Note that
$G$ acts transitively on $\cS$, and that $n:=\#\cS$ is the degree of $f$.

We first give the standard (transparent) group-theoretic translation of
exceptionality.

\begin{lemma}\label{exc gt}
  $f$ is exceptional if and only if the only $A$-orbit on $\cS\times \cS$
  which is also a $G$-orbit is the diagonal.
\end{lemma}

The next two lemmas provide a group-theoretic counting argument which
we use to relate exceptionality to injectivity and surjectivity.
These are variants of \cite[Lemma 6]{Co} and \cite[Lemma 13.1]{FGS}.

\begin{lemma}\label{fixed points}
  Let $A$ be a finite group acting on a finite set $\cT$, let $G$ be a
  normal subgroup of $A$ with $A/G$ cyclic, and let $aG$ be a
  generator of $A/G$.  Then the number of $A$-orbits on $\cT$ which are
  also $G$-orbits equals
\begin{equation*}
\frac{1}{\#G}\sum_{\alpha\in aG}\#\cT^\alpha,
\end{equation*}
where $\cT^\alpha$ denotes the set of fixed points of $\alpha$ on $\cT$.
\end{lemma} 

\begin{proof}
  By examining the different $A$-orbits separately, we may assume that
  $A$ is transitive on $\cT$.  If $G$ is transitive then $ag$ has a
  fixed point for some $g\in G$ (since, if $a$ maps $s\mapsto t$, we
  can choose $g$ mapping $t\mapsto s$); conversely, if some $ag$ has a
  fixed point then the $G$-orbit containing this fixed point must also
  be an $A$-orbit, hence (since $A$ is transitive) must equal $\cT$.
  Thus $G$ fails to be transitive if and only if both sides of the
  equation are zero.  So assume that $A$ and $G$ are both transitive
  (so they have precisely one common orbit).
  
Put
\begin{equation*}
\cV=\{(\alpha,t)\in aG\times \cT\colon \alpha(t)=t\}.
\end{equation*}
Let $A_t$ be the stabilizer of $t$ in $A$ (and similarly for $G$).  On
the one hand, if $\alpha\in aG$ fixes $t$, then $A_t\cap aG=\alpha
G_t$, so there are $\#G_t$ elements in $A_t\cap aG$; hence
$\#\cV=(\#G_t)(\#\cT)=\#G$.  On the other hand, $\#\cV=\sum_{\alpha\in
  aG} \#\cT^\alpha$, and the result follows.
\end{proof}

\begin{lemma}
\label{group theory}
Let $A$ be a finite group acting on a finite set $\cS$, and let $G$ be a
transitive normal subgroup of $A$ with $A/G$ cyclic.  Then the
following properties are equivalent:
\begin{enumerate}
\item The only $A$-orbit on $\cS\times \cS$ which is also a $G$-orbit is
  the diagonal.
\item Every $a\in A$ with $\langle aG\rangle=A/G$ has a unique fixed
  point in~$\cS$.
\item Every $a\in A$ with $\langle aG\rangle=A/G$ has at most one
  fixed point in~$\cS$.
\item Every $a\in A$ with $\langle aG\rangle=A/G$ has at least one
  fixed point in~$\cS$.
\end{enumerate}
\end{lemma}

\begin{proof}
  Since $G$ is transitive on $\cS$, by applying the previous lemma to
  $\cT=\cS$ we see that the average of $\#\cS^\alpha$ (over all $\alpha$ in
  a generating coset of $A/G$) is 1.  Thus (2), (3), and (4) are all
  equivalent to one another.  Applying the previous lemma to
  $\cT=\cS\times \cS$ shows that (1) is equivalent to the average of
  $(\#\cS^\alpha)^2$ being 1.  Since $(\#\cS^\alpha)^2\ge\#\cS^\alpha$, with
  equality if and only if $\#\cS^\alpha$ is 0 or 1, it follows that (1)
  is equivalent to having every $\#\cS^\alpha\le 1$.  Hence (1) and (3)
  are equivalent, which completes the proof.
\end{proof}

Combining these three lemmas with Lemma \ref{from-Wa} yields a proof
that exceptionality implies bijectivity.

\begin{prop}[Lenstra]\label{exc bij}
  If $f$ is exceptional then $f$ is bijective on rational points.
\end{prop}
\begin{proof}
  Let $P \in Y(\Fq)$.  By the definition of finiteness (\cite[p.\ 84]{H}),
  there is an affine open subset $M = \Spec B$ of $Y$ with $P
  \in M(\Fq)$ such that $f^{-1}(M)$ is affine and can be written as
  $\Spec U$ for a ring $U$ that is finite (and therefore integral)
  over $B$.  Since $X$ is normal, $U$ must be the integral closure of
  $B$ in $\Fq(X)$.  Let $R$ be the integral closure of $B$ in
  $\Omega$.  As at the beginning of the section, we let $A =
  \Gal(\Omega / \Fq(Y))$, let $H=\Gal(\Omega/\Fq(X))$, and let $\cS$ be
  the set of left cosets of $H$ in $A$.  Then $R^H = U$ and $R^A = B$.
  Let $\cP$ be the prime in $B$
  corresponding to $P$ and let $D$ and $I$ be the decomposition and
  inertia groups at some prime $\cQ$ of $R$ lying over $\cP$.  Then
  $D/I$ is cyclic and (since $\cP$ has degree one) $DG=A$.  Since $f$ is
  exceptional, Lemmas \ref{exc gt} and \ref{group theory} show that
  every $a\in A$ with $\langle aG \rangle =A$ has a unique fixed point
  in $\cS$.  Since $A=DG$ and $I \subseteq D \cap G$, every $d\in D$
  with $\langle dI \rangle =D$ also satisfies $\langle dG \rangle =A$
  and hence has a unique fixed point in $\cS$.  By Lemma \ref{fixed
    points}, $D$ and $I$ have a unique common orbit on $\cS$.  Thus,
  Lemma \ref{from-Wa} implies there is exactly one maximal ideal
  $\cJ$ in $U$ lying over $\cP$ such that $m_{\cJ} =
  m_{\cP}$, which means there is a unique point $J \in X(\Fq)$
  such that $f(J) = P$.  Hence, $f$ is bijective on rational points.
\end{proof}


We now restrict to the case $\dim X=1$, and prove the converse to
Proposition \ref{exc bij} for sufficiently large $q$.  To give an
explicit bound on $q$ we need the following estimate on the genus
of the Galois closure $\Omega$ of $\Fq(X)/\Fq(Y)$.  Here $n=\deg f$.

\begin{lem}\label{genus bound}
The genus of $\Omega$ satisfies
$$
g_\Omega \le 1+\#G\cdot\frac{g_X-1-(n-2)(g_Y-1)}2 \le 1+
n!\cdot\frac{g_X+n-3}2.$$
\end{lem}
\begin{proof}
Combining Riemann-Hurwitz with Hilbert's formula for the degree of
the ramification divisor gives the formula
$$ 2g_X-2 = n(2g_Y-2) + \sum_{y\in Y} \sum_{i\ge 0} \frac
    {n-\#(\cS\backslash G_i(y))}{[G_0(y):G_i(y)]}, $$
where $G_i(y)$ denotes the $i$-th higher ramification group (in the
lower numbering) of $\Omega/\Fq(Y)$ at a point of $\Omega$ lying over $y$,
and $\cS\backslash U$ denotes the set of orbits of $U$ on $\cS$.  But any
group $U$ of permutations on $\cS$ trivially satisfies
$(\#U)(n-(\#\cS\backslash U)) \ge 2(\#U-1)$, so we get
\begin{align*}
\frac{g_\Omega-1}{\#G} &= g_Y-1 + \frac{1}{2}\sum_{y \in Y}\sum_{i\ge 0}
                          \frac{\#G_i(y)-1}{\#G_0(y)} \\
&\le g_Y-1 + \frac{1}{4}\sum_{y \in Y}\sum_{i\ge 0}
            \frac{n-\#(\cS\backslash G_i(y))}{[G_0(y):G_i(y)]} \\
&= g_Y-1 + \frac{g_X-1}2 - n\frac{g_Y-1}2.
\end{align*}
The result follows.
\end{proof}

\begin{rmk}
The above lemma bounds $g_\Omega$ in terms of $\#G$, and then applies
the trivial bound $\#G\le n!$.  
However, it is often the case that $\#G$ is much smaller than $n!$.
For instance, if $G$ is primitive (meaning that, geometrically, the cover
$X\to Y$ does not have a proper subcover) and $G$ is not $S_n$ or $A_n$,
then $\#G$ is vastly less than $n!$.  Indeed, with explict exceptions,
the order of the group will be polynomial in $n$.  Without exceptions,
one knows that a primitive subgroup of $S_n$ which doesn't contain $A_n$
must have order less than $4^n$ \cite{PS}.

Note also that if $g_X > 1$ then
we have the lower bound
$$
g_{\Omega} \ge 1 +  {\#G}\cdot (g_X -1)/n,
$$
so our upper bound has the right order of magnitude in this situation.
\end{rmk}


Combining the above lemmas with Corollary~\ref{Lenstra Cheb} yields the
main result of this section:
 
\begin{thm}\label{part}
  Let $f: X \lra Y$ be a finite, separable morphism of smooth projective
  curves and let $n$ be the degree of $f$.  Suppose that $\sqrt{q} \ge
  n!(3g_X + 3 n)$ and that $f$ is either injective or surjective on
  rational points.  Then $f$ is exceptional, and is bijective on
  rational points.
\end{thm}

\begin{proof}
  Let $a$ be an element of $A$ such that $\langle aG\rangle=A/G$.
  Let $\cU$ be the set of points in $Y(\Fq)$ over which
  $f$ ramifies.  By Riemann-Hurwitz, $\# \cU \leq (2g_X + 2n - 2)$.
  Combining this bound with Lemma \ref{genus bound} and the inequality
  $[\Omega:\Fq(Y)] \leq n!$, we see that
\begin{equation*}
   2 g_\Omega + \sqrt{[\Omega:\Fq(Y)] (\#\cU)} < n!(3 g_X + 3n).
  \end{equation*}
  Now Corollary~\ref{Lenstra Cheb} implies there is a closed point $Q$
  of $\Omega$ which lies over a point $P\in Y(\Fq)$ such that $Q/P$ is
  unramified and its decomposition group is generated by $a$.  By
  Lemma~\ref{from-Wa}, the number of points in $X(\Fq)$ lying over $P$
  equals the number of points of $\cS$ fixed by $a$.
  
  Thus, surjectivity of $f$ implies property (4) of Lemma~\ref{group
    theory}, and injectivity of $f$ implies property (3).  By
  Lemma~\ref{group theory} and Lemma~\ref{exc gt}, if $f$ is either
  surjective or injective then $f$ is exceptional, and hence (by
  Proposition~\ref{exc bij}) $f$ is bijective.
\end{proof}

\begin{rmk}
  The above result (and its proof) remains valid under the weaker
  hypothesis that the map $X(\Fq)\to Y(\Fq)$ is either injective or
  surjective over non-branch points of $Y$.  More generally, let $\cU$
  be any subset of $Y(\Fq)$ which includes all points of $Y(\Fq)$ over
  which $f$ is ramified.  The above proof shows that either
  injectivity or surjectivity of $X(\Fq)\setminus f^{-1}(\cU)\to
  Y(\Fq)\setminus\cU$ implies exceptionality of $f$ (and hence
  bijectivity of $X(\Fq)\to Y(\Fq)$), so long as $\sqrt{q}\ge
  2g_\Omega+\sqrt{\#G\cdot\#\cU}$.
\end{rmk}


\section{Examples and Further Directions}\label{examples}

We first give examples of covers of curves $X\to Y$ over $\Fq$ which
are surjective but not injective on rational points, as well as
examples which are injective but not surjective.  In these examples,
the degree $n$ of the cover is small relative to $q$, but the genus of
$X$ is as large as $q$.
\begin{eg}
  Let $q$ be an odd prime power and let $n > 1$ divide $(q-1)/2$.
  Pick $a,\gamma\in\Fq^*$ with $a$ an $n$-th power.  Let $X$ be the
  normalization of the affine curve
  $$
  y^n = \gamma \prod\limits_{t \in \Fq^* \setminus a}(x-t),$$
  so
  $X$ has genus $(n-1)(q-3)/2$.  Let $f:X \lra \bP^1_{\Fq}$ be the
  morphism induced by projection onto the $x$-coordinate.  Then $f$ is
  totally ramified over all the rational points of $\bP^1_{\Fq}$
  (including infinity) except for 0 and $a$.  Moreover, $\prod_{t\in
    \Fq^* \setminus a}(0-t)=a^{-1}$ and $\prod_{t\in\Fq^* \setminus
    a}(a-t)=-a^{-1}$ are both $n$-th powers in $\Fq^*$.  Thus, if
  $\gamma$ is an $n$-th power in $\Fq^*$, then both $x=0$ and $x=a$
  split completely under $f$, so $f$ is surjective on rational points
  but not injective.  If $\gamma$ is not an $n$-th power, then neither
  $x=0$ nor $x=a$ is the image of an $\Fq$-point under $f$, so $f$ is
  injective on rational points but not surjective.\qed
\end{eg}

Our next two examples are functions of the same shape, all of which
are bijective but only some of which are exceptional.  The exceptional
ones come from the largest known family of exceptional functions.

\begin{eg}
  We consider the case of degree-5 maps between genus-0 curves.  It
  follows from Dickson's results~\cite{Di} that, if such a map is
  totally ramified over some $\Fq$-rational point and $q>13$, then
  injectivity on rational points implies exceptionality.  Here are
  examples showing the necessity of the ramification hypothesis: the
  rational function $(x^5-ax)/(x^4-b)$ is non-exceptional and
  bijective on $\Fq$-points if $(q,a,b)$ is either $(17,10,3)$ or
  $(29,13,4)$.\qed
\end{eg}

\begin{eg}
  The above example has the same shape as some exceptional maps.
  Namely, if $k$ is any field containing a primitive fourth root of
  unity $i$ and a nonsquare $b$, then $f(x)=(x^5-b(4i-3)x)/(x^4-b)$ is
  an exceptional map $\bP^1\to\bP^1$ over $k$.  These examples come
  from the construction in \cite{Fr2} as follows.  Let $\varphi\colon
  E\to E'$ be the 5-isogeny between the elliptic curves $E\colon
  w^2=x^3+xb(1+2i)/4$ and $E'\colon v^2=u^3+ub(1+2i)^5/4$ such that
  the nontrivial elements in the kernel of $\varphi$ are the pairs
  $(x,w)$ with $x^2=-b/4$ and $w^2=xbi/2$.  Map $E$ and $E'$ to
  $\bP^1$ by taking the quotient by the automorphism (of curves)
  $P\mapsto (0,0)-P$.  Then $\varphi$ induces a map $\bP^1\to\bP^1$ which
  is easily seen to be our $f$.  More generally, the largest known
  supply of exceptional maps $\bP^1\to\bP^1$ are maps induced from
  isogenies of elliptic curves (cf.\ \cite{Fr2,GMS}).\qed
\end{eg}

We next give a large class of exceptional covers of curves.

\begin{eg}
  Let $C$ be a curve on an abelian variety $A$ over $\Fq$, let
  $\varphi\colon A\to A$ be the multiplication-by-$d$ map, and suppose
  $d$ is coprime to both $q$ and $\#A(\Fq)$.  Suppose furthermore that
  $\varphi^{-1}(C)$ is geometrically irreducible.  Then the map
  $\varphi^{-1}(C)\to C$ is exceptional.  This follows from the fact
  that the induced map from $\varphi^{-1}(C)$ to $C$ is bijective over
  any extension of $\ell$ of $\Fq$ that does not contain the field of
  definition of any of the points in $A({\overline \F}_q)$ having order
  a nontrivial divisor of $d$.
\end{eg}

We conclude by discussing possible higher-dimensional analogs of
Theorem \ref{the main}.

\begin{conj}\label{main-conj}
  Let $f\colon X \lra Y$ be a finite, generically \'{e}tale map of
  degree $n \geq 2$ between two smooth projective varieties of
  dimension $r$ defined over $\Fq$.  Then there exists a constant $C$,
  depending only on $n$, the dimensions of $X$ and $Y$, and the
  $\ell$-adic Betti numbers $b_1,\dots,b_{2r-1}$ of $X$, such that if
  $q > C$ and $f$ induces an injection or a surjection from $X(\Fq)$
  to $Y(\Fq)$, then $f$ is exceptional and gives a bijection from
  $X(\Fq)$ to $Y(\Fq)$.
\end{conj}

We can prove Conjecture~\ref{main-conj} for maps $f: \bP^m \lra \bP^m$
(we mean the implication of exceptionality since injectivity and
surjectivity are equivalent) although we are not able to give a simple
formula for $C$.  Here is a sketch of the proof.  If $D$ is a
geometric component of $\bP^m \times_f \bP^m$ defined over $\Fq$, then
$D$ is birational to a subvariety of $\bP^{2m}$ of dimension $m$ and
degree at most $(2n)^m$.  Thus, by Lang-Weil \cite{LW},
there is a constant $C_1$, depending
only on $n$ and $m$, such that $\#D(\Fq) \geq q^m - C_1 q^{m - 1/2}$.
Arguing as in Proposition~\ref{INT}, we also see that, since $f$ is
injective, $\#D(\Fq) \leq R(\Fq)$ where $R$ is the ramification locus
of $f$.  Since $R$ is a divisor of degree at most $2 m n$ on $\bP^m$,
there is a constant $C_2$, depending only on $m$ and $n$ such that
$\#R(\Fq) \leq C_2 q^{m-1}$ (again by Lang-Weil), which contradicts our
earlier lower bound on $\#D(\Fq)$ when $q$ is sufficiently large.
  
Unfortunately, proving Conjecture~\ref{main-conj} in general seems to
be much more complicated since we cannot use Lang-Weil and are
forced instead to attempt to control Betti numbers of various
varieties that arise.  One possibility, suggested by
Lucien Szpiro, is to directly prove the equivalence of
injectivity and surjectivity in higher dimensions by examining the
induced maps from curves in $X$ to curves in $Y$.

The best we can do for maps between general varieties is the following
non-explicit version of Conjecture~\ref{main-conj}, where we allow
the constant $C$ to depend on the map $f$.

\begin{prop}\label{end}
  Let $f:X \lra Y$ be a finite separable map between normal
  varieties over $\Fq$.  If $f$ induces an surjective or injective map
  from $X(\Fqm)$ to $Y(\Fqm)$ for infinitely many $m$, then $f$ is
  exceptional.
\end{prop}

\begin{proof}
  Let $D_f$ denote the ramification locus of $f$ in $X$.  Let $W = X
  \setminus D_f$, and let $V$ be the normalization of $W$ in $\Omega$
  (the Galois closure of $\Fq(X)$ over $\Fq(Y))$).  Let $A_m =
  \Gal(\Fqm(\Omega) / \Fqm(W))$ and let $G_m =\Gal(\ell(\Omega) /
  \ell(W))$, where $\ell$ is the closure of $\Fqm$ in $\Omega.\Fqm$.
  By Corollary~\ref{non smooth}, there exists $M$ such that for any $m
  \geq M$ and any $\sigma \in A_m$ such that $\langle \sigma G_m
  \rangle = A_m/G_m$, there is a point $P \in W(\Fqm)$ and a closed
  point $Q$ on $V_{\Fqm}$ such that $\langle \sigma I(Q/P) \rangle =
  D(Q/P)/I(Q/P)$.  For such $m$, injectivity and surjectivity each
  imply exceptionality by Lemma~\ref{group theory}.
\end{proof} 
   
We note that Fried (\cite{Fr}) has previously proved
Proposition~\ref{end} above in the case of maps from affine space
to itself.



\begin{thebibliography}{99}
\newcommand{\au}[1]{{#1},}
\newcommand{\ti}[1]{\textit{#1},}
\newcommand{\jo}[1]{{#1}}
\newcommand{\vo}[1]{\textbf{#1}}
\newcommand{\yr}[1]{(#1),}
\newcommand{\pp}[1]{#1.}
\newcommand{\pps}[1]{#1;}
\newcommand{\bk}[1]{{#1},}
\newcommand{\inbk}[1]{in: {#1},}
\newcommand{\xxx}[1]{{arXiv:#1}}

\bibitem[Bo]{Bo} 
\au{N. Bourbaki}
\bk{Elements of Mathematics, Commutative Algebra}
Hermann, Paris, France, 1972.
 
\bibitem[Co]{Co}
\au{S. D. Cohen}
\ti{The distribution of polynomials over finite fields}
\jo{Acta Arith.}
\vo{17}
\yr{1970}
\pp{255--271}

\bibitem[DL]{DL}
\au{H. Davenport and D. J. Lewis}
\ti{Notes on congruences. I}
\jo{Quart. J. Math. Oxford}
\vo{14}
\yr{1963}
\pp{51--60}

\bibitem[De]{Deligne} 
\au{P. Deligne}
\ti{La conjecture de Weil. I}
Inst. Hautes \'{E}tudes Sci. Publ. Math. No. 43 (1974), 273--307. 

\bibitem[Di]{Di}
\au{L. E. Dickson}
\ti{The analytic representation of substitutions on a power of a prime
  number of letters, with a discussion of the linear group}
\jo{Ann. of Math.}
\vo{11}
\yr{1896}
\pp{65--120}

\bibitem[Fr]{Fr}
\au{M. Fried}
\ti{On a theorem of MacCluer}
\jo{Acta Arith.}
\vo{25}
\yr{1974}
\pp{121--126}

\bibitem[Fr2]{Fr2}
\au{M. Fried}
\ti{Galois groups and complex multiplication}
\jo{Trans. Amer. Math. Soc.}
\vo{235}
\yr{1978}
\pp{141--163}

\bibitem[Fr3]{Fr3}
\au{M. Fried}
\ti{Global construction of general exceptional covers}
\inbk{Finite Fields: Theory, Applications, and Algorithms}
American Mathematical Society, Providence
\yr{1994}
\pp{69--100}

\bibitem[FGS]{FGS}
\au{M. Fried, R. Guralnick and J. Saxl}
\ti{Schur covers and Carlitz's conjecture}
\jo{Israel J. Math.}
\vo{82}
\yr{1993}
\pp{157--225}

\bibitem[FJ]{FJ}
\au{M. Fried and M. Jarden}
\bk{Field Arithmetic}
Springer-Verlag,
Berlin, 1986.

\bibitem[GM]{GM}
\au{J. von zur Gathen and K. Ma}
\ti{The computational complexity of recognizing permutation functions}
\jo{Comput. Complexity}
\vo{5}
\yr{1995}
\pp{76--97}

\bibitem[Gr]{EGA}
\au{A. Grothendieck}
\ti{\'{E}l\'{e}ments de g\'{e}om\'{e}trie alg\'{e}brique. IV. \'{E}tude locale
  des sch\'{e}mas et des morphismes de sch\'{e}mas. II}
Inst. Hautes \'{E}tudes Sci. Publ. Math. No. 24 (1967), 1--231.

\bibitem[GMS]{GMS}
\au{R. Guralnick, P. Mueller and J. Saxl}
\ti{The rational function analogue of a question of Schur and exceptionality
of permutation representations}
\jo{Mem. Amer. Math. Soc.}
\vo{162},
no. 773 (2003), 1--79.

\bibitem[Ha]{H}
\au{R. Hartshorne}
\bk{Algebraic Geometry}
Springer-Verlag,
New York, 1977.

\bibitem[LW]{LW}
\au{S. Lang and A. Weil}
\ti{Number of points on varieties in finite fields}
\jo{Amer. J. Math.}
\vo{76}
\yr{1954}
\pp{819--827} 


\bibitem[LY]{LY}
\au{D. Leep and C. Yeomans}
\ti{The number of points on a singular curve over a finite field}
\jo{Arch. Math.}
\vo{63}
\yr{1994}
\pp{420--426}

\bibitem[LMZ]{LMZ}
\au{H. W. Lenstra, Jr., D. Moulton, M. Zieve}
\ti{Exceptional covers}
in preparation.

\bibitem[Mac]{Ma}
\au{C. R. MacCluer}
\ti{On a conjecture of Davenport and Lewis concerning exceptional polynomials}
\jo{Acta Arith.}
\vo{12}
\yr{1967}
\pp{289--299}

\bibitem[Mat]{Mats}
\au{H. Matsumura}
\bk{Commutative ring theory}
Cambridge University Press, Cambridge, 1986.

\bibitem[MS]{MS}
\au{V. K. Murty and J. Scherk}
\ti{Effective versions of the Chebotarev density theorem for function fields}
\jo{C.R. Acad. Sci. Paris S\'{e}r. I Math.}
\vo{319}
\yr{1994}
\pp{523--528}

\bibitem[PS]{PS}
\au {C. Praeger and J. Saxl}
\ti {On the orders of primitive permutation groups}
\jo{Bull. London Math. Soc.}
\vo{12} 
\yr{1980} 
\pp{303--307}


\bibitem[Wa]{Wa}
\au{B. L. van der Waerden}
\ti{Die Zerlegungs- und Tr\"agheitsgruppe als Permutationsgruppen}
\jo{Math. Ann.}
\vo{111}
\yr{1935}
\pp{731--733}

\bibitem[We]{We}
\au{A. Weil}
\bk{Vari\'et\'es ab\'eliennes et courbes alg\'ebriques}
Hermann, Paris, 1948.

\bibitem[Wi]{Wi}
\au{K. S. Williams}
\ti{On exceptional polynomials}
\jo{Canad. Math. Bull.}
\vo{11}
\yr{1968}
\pp{279--282}

\end{thebibliography}
\end{document}